\documentclass[a4paper, 10pt, draft, reqno]{amsart}
\usepackage{epsf}
\usepackage{amsmath}
\usepackage{amssymb}
\usepackage{latexsym}

\usepackage[cmtip,matrix,arrow]{xy}
\UseComputerModernTips
\CompileMatrices

\DeclareMathOperator{\sgn}{sgn}
\newcommand{\Z}{{\mathbb Z}}

\newcommand{\N}{{\mathbb N}}
\newcommand{\DD}{{\mathcal D}}

\DeclareMathOperator{\Hom}{Hom}
\DeclareMathOperator{\Ext}{Ext}

\DeclareMathOperator{\rad}{rad}

\newcommand{\frob}{^{\mbox{\rm\tiny F}}}

\newcommand{\otherwise}{\mbox{\rm otherwise}}
\newcommand{\wif}{\mbox{\rm if }}
\newcommand{\wand}{\mbox{\rm and }}

\newcommand{\Mod}{\mathrm{mod}}
\newcommand{\MMod}{\mathrm{Mod}}

\newcommand{\GL}{\mathrm{GL}}
\newcommand{\SL}{\mathrm{SL}}

\newcommand{\Sym}{\mathfrak{S}}
\DeclareMathOperator{\partn}{\Lambda^+}
\DeclareMathOperator{\cpartn}{\Lambda^+_{\mbox{\scriptsize{\rm{col}}}}}
\DeclareMathOperator{\rpartn}{\Lambda^+_{\mbox{\scriptsize{\rm{row}}}}}

\newcommand{\bnabla}{\overline{\nabla}}

\begin{document}
\theoremstyle{plain}
\newtheorem{thm}{Theorem}[section]
\newtheorem{propn}[thm]{Proposition}
\newtheorem{lem}[thm]{Lemma}
\newtheorem{cor}[thm]{Corollary}
\newtheorem{conj}[thm]{Conjecture}
\theoremstyle{definition}
\newtheorem{rem}[thm]{Remark}
\newtheorem{ass}[thm]{Assumption}
\newtheorem{defn}[thm]{Definition}
\newtheorem{example}[thm]{Example}

\setlength{\parskip}{1ex}

\title{Homomorphisms and Higher Extensions
for Schur algebras and symmetric groups}
\author{Anton Cox
\and Alison Parker} 
\address{Centre for Mathematical Science \\
City University  \\
Northampton Square, London, EC1V 0HB\\
 England.} 
\email{A.G.Cox@city.ac.uk}
\address{
School of Mathematics and Statistics F07\\
University of Sydney \\
NSW 2006\\
Australia.}
\email{alisonp@maths.usyd.edu.au}

\begin{abstract}
This paper surveys, and in some cases generalises, many of the recent
results on homomorphisms and the higher $\mathrm{Ext}$ groups for
$q$-Schur algebras and for the Hecke algebra of type $A$.  We review
various results giving isomorphisms between Ext groups in the two
categories, and discuss those cases where 
explicit results have been determined.  
\end{abstract}
\maketitle

\section{Introduction}

Since 1901 and the work of Schur \cite{schur} it has been known that
the representation theories of the symmetric and general linear groups
are intimately related. Following later work of Schur \cite{schurb}
and Carter and Lusztig \cite{carlus}, Green \cite{green} showed how
this relationship could be explained via passage through the Schur
algebra using the Schur functor. This algebra has been an object of intense
study ever since, and with the rise of quantum groups was generalised
by Dipper and James \cite{dj1} to give the $q$-Schur algebra. This
plays a corresponding role in relating the representation theories of
type $A$ Hecke algebras and the quantum general linear group.

On the level of cohomology however, there are  some striking
differences between the symmetric groups and the Schur algebras. For
example, Schur algebras are quasi-hereditary and hence have finite
global dimension, while the global dimension of the symmetric group is
usually infinite. Thus we must expect some difficulties to arise in
the passage of cohomological data from one category to the other. 

The majority of this paper is concerned with the $\Ext$-groups
corresponding to pairs of induced (or Weyl) modules for the
algebraic/quantum group, and to pairs of Specht modules for the
symmetric group/Hecke algebra. In their respective categories these
classes of modules play a key role, and correspond under the action of
the Schur functor, so we might hope that their cohomology theories
are closely related. We will review what is known, and show how some
of the results can be extended and refined.

Carter and Lusztig \cite{carlus} showed that at the level of
homomorphisms these two theories could essentially be identified under
the standard correspondence. Surprisingly however, consideration of
the relationship between the higher $\Ext$-groups has been a relatively
recent phenomenon.

The first real results in that direction are due to Doty, Erdmann and
Nakano \cite{den}, who used certain spectral sequences arising in
the general setting of finite dimensional algebras with idempotent
functors. These gave rise to various relationships between certain
$\Ext^1$ and $\Ext^2$ spaces. Using this framework, Kleshchev and
Nakano \cite{klenak} and Hemmer and Nakano \cite{hemnak} were able to
refine these results further, and recently Donkin \cite{donkpreprint}
and Parshall and Scott \cite{psqweyl} have generalised these results. 

Donkin, and Parshall and Scott, show that there are equalities of
$\Ext$-groups in degrees greater than zero (for certain classes of
modules) provided that the degree is small compared with the
characteristic of the field (or the degree of the root of unity in the
quantum case). This will be our key tool in relating the two
categories. 

There is an important reduction result due to Donkin
\cite{donklnm} (surprisingly little-known), which can be used
to reduce to lower rank calculations in certain special cases, and a
similar result due to Cline, Parshall, and Scott \cite{cpstransf}
where we allow one of the modules to be simple.

Unfortunately, relatively little is known about the actual cohomology
for explicit modules in either category. We do have an upper bound on
the degree of a non-zero $\Ext$ group between two non-zero modules
given by \ref{ryom} and \ref{parker}, but this is the best we can do
in general. 
On the algebraic group side
the only setting where complete information is known is $\SL_2$ 
(and $q$-$\GL_2(k)$), where
the second author \cite{par3} has recently determined all $\Ext$-groups
between Weyl modules. This extends work of Erdmann and the first
author \cite{erd, cox, coxerd} where the groups up to degree $2$ where
determined. For $\SL_3$, the authors have determined all homomorphisms
between Weyl modules in characteristic at least $3$ \cite{coxpar}.
(This has recently been generalised to $q$-$\GL_3(k)$ in \cite{par4}.)

Apart from these low-rank calculations, the other main general results
are either for weights which are close together (for example in
\cite{andGBcoh}, \cite{kop1} and \cite{wen}) or for weights which differ by a
single reflection. In the latter case there is a theorem of
Carter and Payne \cite{cpmaps} which has recently been generalised by
Fayers and Martin \cite{faymar}. This generalisation gives a set of
sufficient conditions for the existence of a homomorphism which can be
difficult to verify in particular examples, and we will show how these
conditions may be simplified.

Other results in this area consider extensions between simple modules.
Kleschev and Sheth \cite{klsh} consider $\Ext^1$ between simple
modules for the symmetric group. They show  in some cases that
an extension between two simple modules $D^\lambda$ and $D^\mu$
must appear in the Specht module $S^\mu$ when $\lambda$ and $\mu$ have
at most $p-1$ parts and $p \ge 3$. We will give a generalisation of
this result. 

Wherever possible we have stated results using partition combinatorics
rather than alcove geometry and the affine Weyl group, in order to
simplify the exposition. However, we believe that the alcove-geometric
approach has been under-utilised in the Hecke algebra setting and
provides a natural context for many results in this area (see for
example Section \ref{homsec} and the references therein).

\section{Notation}\label{prelim}

Let $k$ be an algebraically closed field of characteristic $p\geq 0$
and $q\in k$ a root of unity.  We define $l$ to be $p$ if $q=1$, or
the smallest non-zero exponent for which $q^l =1$ if $q \ne 1$.  We
are interested in the relationship betwen the representation theories
of two objects: $G=q$-$\GL_n(k)$, the quantum general linear group as
defined by Dipper and Donkin \cite{dipdonk}, and $H=H_{q,n}(k)$, the
Hecke algebra of type $A$.  If $q=1$ then $G$ is the classical general
linear group $\GL_n(k)$, the group of $n$ by $n$ invertible matrices
over $k$, and $H$ is the group algebra of the symmetric group
$\Sym_r$. In this case we will further assume that $p>0$.  For the
rest of this section we will review the basic results which we
require.  This material can be found in \cite{donkbk} for the quantum
group and \cite{mathasbk} for the Hecke algebra (which are
respectively in the spirit of the classical expositions in
\cite{jantzen} and \cite{jk}).

Our main tool for relating the two module categories will be the
$q$-Schur algebra $S(n,r)$, as defined in \cite{dj1} (generalising the
classical definition of Green \cite{green}). This will often be
abbreviated to $S$ when the context is clear.  The module category of
$S(n,r)$ can be identified with the category of polynomial modules for
$G$ of degree $r$. If $r\leq n$ there exists a certain idempotent
$e\in S$. This give rise to the Schur functor $f:S\longrightarrow eSe$
defined on modules by $f(M)=eM$, together with a partial inverse $g$
given by $g(N)=Se\otimes_{eSe} N$. As $eSe\cong H$ this provides the
desired connection between the representations of $G$ and $H$.

We first consider the quantum group. Let $X \cong \Z^n$ be the set of
weights for $G$. We chose a Borel $B$ such that the dominant weights
are given by $X^+ = \{ (\lambda_1, \lambda_2, \ldots, \lambda _n) \in
X \mid \lambda_1 \ge \lambda _2\ge \cdots \ge \lambda_n\}$. When a
dominant weight $\lambda$ satisfies $\lambda_n\geq 0$ we call it a
\emph{polynomial weight}, and identify it with the corresponding
partition of $|\lambda| = \sum_{i} \lambda _i$. We denote by
$\MMod(G)$ the rational modules for $G$ and by $\Mod(G) \subset
\MMod(G)$ the finite dimensional rational modules for $G$.  We have
two dualities on $\Mod(G)$; the usual linear dual $^*$, and a
contravariant duality $^{\circ}$ (which in the classical case is
induced by the transpose map on matrices).

Given a dominant weight $\lambda$ there exists a corresponding
one-dimensional $B$-module $k_{\lambda}$. We let $\nabla(\lambda)$
denote the $G$-module obtained from this by induction. This has
highest weight $\lambda$, and simple socle which we denote by
$L(\lambda)$. These socles give a full set of inequivalent simple
$G$-modules.  The contravariant dual
$\Delta(\lambda)=\nabla(\lambda)^{\circ}$ is called a \emph{Weyl module}.

Let $\partn(r)$ denote the set of partitions of $r$, and $\partn(n,r)$
denote the subset of partitions of $r$ with at most $n$ parts. Given a
partition $\lambda$ we denote the conjugate partition by $\lambda'$. A
partition $\lambda$ of $r$ is \emph{column $l$-regular} if $\lambda_i
- \lambda_{i+1} < l$ for all $i$. A partition $\lambda$ of $r$ is
\emph{row $l$-regular} if and only if $\lambda'$ is column
$l$-regular. We let $\cpartn(n,r)$ (respectively $\rpartn(n,r)$)
denote the subset of $\partn(n,r)$ consisting of the column
(respectively row) $l$-regular partitions. Let $\mu$ and $\lambda$ be
two partitions of $r$.  We say $\lambda$ \emph{dominates} $\mu$ if
$|\lambda|=|\mu|$ and $\sum_{i =1}^{j} \mu_i \le \sum_{i=1}^{j}
\lambda_i $ for all $ 1 \le j \le n$ and write $\mu \le \lambda$.

The simple modules for $S(n,r)$ are precisely the $L(\lambda)$ for $G$
with $\lambda\in \partn(n,r)$. We define a symmetric power associated
to each $\lambda \in \partn(n,r)$ by
$$S^\lambda E = S^{\lambda_1}E \otimes S^{\lambda_2} E \otimes \cdots
\otimes S^{\lambda_n} E$$ where $S^i E$ is the $i$th symmetric power
of the natural module $E$.  The $S^\lambda E$ are injective
$S(n,r)$-modules, but are not indecomposable in general.  We let
$I(\lambda)$ be the injective hull of $L(\lambda)$, regarded as an
$S$-module.  The module $I(\lambda)$ is a direct summand of $S^\lambda
E$ and is the unique summand of $S^\lambda E$ with $L(\lambda)$ as
its socle. 

Although we have a labelling of the simples, they are not well
understood.  We do not even know their dimension in general.  To
determine their characters it would be enough to determine
$[\nabla(\lambda):L(\mu)]$ for all pairs $\mu\leq \lambda$, where
$[M:L]$ denotes the multiplicity of the simple module $L$
as a composition factor of $M$.

We next consider the Hecke algebra $H_{q,n}(k)$. For each
$\lambda\in\partn(r)$ we can define corresponding $H$-modules
$S^{\lambda}$, the \emph{Specht module}, $Y^{\lambda}$, the
\emph{Young module}, and $M^{\lambda}$, the \emph{permutation module}. If
$\lambda \in \rpartn(n,r)$ then $S^{\lambda}$ has simple head, which
we denote by $D^{\lambda}$. These simple heads provide a full set of
inequivalent $H$-modules.

For our purposes these various modules can be defined in terms of the
Schur functor $f$. In particular, if $r\leq n$ we have 
\begin{align}\label{fdoes}
f S^{\lambda}\! E \cong M^{\lambda} \quad\quad
f I(\lambda) \cong Y^{\lambda} \quad\quad
f \nabla(\lambda) \cong S^{\lambda}
\end{align}

We also know that $fL(\lambda)=0$ unless $\lambda\in\cpartn(n,r)$;
however in order to describe the effect of $f$ on the remaining simple
modules we need a little more notation.
We write $\sgn$ for the $q$-analogue of the one-dimensional sign
representation of $\Sym_r$. Then the module $D^\mu \otimes \sgn$ is
again a simple $H$-module, denoted $D^{m(\mu)}$.  The map $m(-)$ on
$\rpartn(n,r)$ is an involutory bijection and is known as the
\emph{Mullineux map}. It can be shown that $f L(\lambda) = D^{m(\lambda')}$
if $\lambda \in \cpartn(n,r)$.

Since $f$ is exact we also have
$[\nabla(\lambda) : L(\mu)] = [S^\lambda :
D^{m(\mu')}]$ for $\mu \in
\cpartn(n,r)$.
We can also show that $Y^\lambda$ is projective  as an
$H$-module
and $I(\lambda)$ is projective as a $S(n,r)$-module if
$\lambda \in
\cpartn(n,r)$.

\section{Equating Ext groups for the Schur algebra and the
  symmetric group}

We wish to investigate the extent to which $\Ext$ groups for our two
categories can be identified. The following fundamental theorem will
be our starting point. It was originally proved in the classical case
by Carter and Lusztig \cite[Theorem 3.7]{carlus} and later generalised
to the quantum case by Dipper and James \cite[8.6
Corollary]{dipjamesq}.

\begin{thm}[Carter-Lusztig, Dipper-James]\label{thm:equalhom}
If $l \ge 3 $ or if $l=2$ and either $\lambda \in \rpartn(n,r)$ or $\mu
\in \cpartn(n,r)$ then
$$\Hom_{S(n,r)}(\nabla(\lambda), \nabla(\mu)) \cong
\Hom_{H}(S^\lambda, S^\mu).$$
\end{thm}

\begin{rem}
This theorem appears in a slightly different form in the literature:
stated in terms of Weyl modules and without the case $\mu
\in \cpartn(n,r)$. The translation from Weyl to induced modules is
clear by applying contravariant duality, while the additional case
follows from the case $\lambda \in \rpartn(n,r)$ and the following
pair of results.
\end{rem}

\begin{lem}[Donkin]\label{lem:Nn}
Suppose that $n \le N$ and $\lambda$, $\mu \in \partn(n,r)$, which is
considered in the usual way as a subset of $\partn(N,r)$.
Then for all $i\geq 0$ we have 
$$\Ext^i_G(\nabla(\lambda), \nabla(\mu)) 
\cong \Ext^i_{q-\GL_N(k)}(\nabla(\lambda), \nabla(\mu)).$$
\end{lem}
\begin{proof}
This follows using \cite[4.2 (17)]{donkbk}.
\end{proof}

Note that it is only the $\Ext$'s
between induced modules that are stable; the extensions between
simples, for example, are dependent on $n$. However if $n$ is larger
than the degree of the partition labelling the simples, then the
$\Ext$ groups become stable (using a result of \cite{green}). 
We will also see that
$\Ext^i_G(L(\lambda), \nabla(\mu)) \cong 
\Ext^i_{q-\GL_N(k)}(L(\lambda), \nabla(\mu))$, this follows
from Theorem \ref{cpscut}.

\begin{lem}[Donkin]\label{conj}
Suppose $r \le n$. Then for all $m \in \N$ we have
$$\Ext^m_{S(n,r)}(\nabla(\lambda), \nabla(\mu))
\cong 
\Ext^m_{S(n,r)}(\nabla(\mu'), \nabla(\lambda'))$$
and
$$\Ext^m_{H}(S^\lambda, S^\mu)
\cong 
\Ext^m_{H}(S^{\mu'}, S^{\lambda'}).$$
\end{lem}
\begin{proof}
The result in the Schur algebra is proved using the Ringel dual,
\cite[Corollaries 3.8 and 3.9]{dontilt}  and
\cite[Proposition 4.1.5]{donkbk}.
If we distinguish the modules for the Ringel dual of $S(n,r)$ by
primes
then we have
\begin{align*}
\Ext^m_{S(n,r)}(\nabla(\lambda), \nabla(\mu))
&\cong 
\Ext^m_{S(n,r)'}(\Delta'(\lambda'), \Delta'(\mu'))\\
&\cong 
\Ext^m_{S(n,r)}(\Delta(\lambda'), \Delta(\mu'))\\
&\cong 
\Ext^m_{S(n,r)}(\nabla(\mu'), \nabla(\lambda'))
\end{align*}
where the first isomorphism follows by \cite[A4.8 (i)]{donkbk}
the second by $S(n,r)$ being Ringel self-dual if $r \le n$
and the third by taking contravariant duals.

The Hecke algebra version follows easily from the 
isomorphism  \cite[Proposition 4.5.9]{donkbk}.
\begin{equation}\label{Sconj}
(S^\lambda)^* \cong
S^{\lambda'}\otimes \sgn.
\end{equation}
We have
\begin{align*}
\Ext^m_H(S^\lambda, S^\mu)
&\cong 
\Ext^m_H((S^\mu)^*, (S^\lambda)^*)\\
&\cong 
\Ext^m_H(S^{\mu'}\otimes \sgn, S^{\lambda'}
\otimes \sgn)\\
&\cong 
\Ext^m_H(S^{\mu'}, S^{\lambda'})
\end{align*}
\end{proof}

We would like to have an analogue of Theorem \ref{thm:equalhom} for
the higher Ext-groups. However, it is easy to see that they cannot
always be equal. We know that $H$ is a symmetric algebra and for
$l\leq r$ is usually not semi-simple. (It is semi-simple if $l > r$.)
In general $H$ has infinite global dimension: i.~e.  for all $m \in
\N$, there exist $M$, $N \in \Mod H$ such that $\Ext^m_{k \Sym_r}(M,N)
\ne 0$.  On the other hand $S(n,r)$ is a quasi-hereditary algebra,
which is usually not symmetric and not semi-simple.  Quasi-heredity
implies that $S(n,r)$ has finite global dimension. So there exists $m
\in \N$ such that $\Ext^j_{S(n,r)}(M,N) = 0$ for all $j > m$ and all $M$, $N
\in \Mod S(n,r)$.

The relationship between the two cohomology theories has been studied
in \cite{den}, in a general setting. Based on this work, Theorem
\ref{thm:equalhom} has been generalised by Kleshchev and Nakano
\cite{klenak} and Hemmer and Nakano \cite{hemnak}, and further
generalised by Donkin \cite[Section 10, Propositions 2 and
3]{donkpreprint} and Parshall and Scott \cite[Theorem
4.6]{psqweyl}. We state here Donkin's version, 
that of Parshall and Scott is similar.

\begin{thm}[Donkin]\label{thm:donk} 
Suppose that $X \in \Mod H$  has a Specht
series 
and that $Y \in \Mod S(n,r)$.
If $l \ge 4 $ and $0 \le i \le l-3$ then 
$$\Ext^i_{S(n,r)}(gX, Y)
\cong \Ext^i_{H}(X, fY).$$
and $g S^\lambda = \nabla(\lambda)$.
\end{thm}

Combining Theorems \ref{thm:equalhom} and \ref{thm:donk} with the
identifications in and after (\ref{fdoes}) we obtain

\begin{cor}[Donkin, Parshall-Scott]\label{cor:sameext}
If $\lambda$, $\mu \in \partn(n,r)$ and $l\ge 3$ then 
$$\Ext^i_{S(n,r)}(\nabla(\lambda), \nabla(\mu))
\cong \Ext^i_{H}(S^\lambda, S^\mu)$$
and
$$\Ext^i_{S(n,r)}(\nabla(\lambda), L(\mu))
\cong \Ext^i_{H}(D^{\mu'},S^{\lambda'})$$
for $0 \le i \le l-3$.
\end{cor}
\begin{proof}
The first isomorphism is immediate from the Theorem; for the second
(which is not explicitly stated by Donkin) we have for $l>3$ that
\begin{align*}
\Ext^i_{S(n,r)}(\nabla(\lambda), L(\mu))& \cong
\Ext^i_H(S^{\lambda}, D^{m(\mu')})\\
& \cong\Ext^i_H(D^{m(\mu')},(S^{\lambda})^*)\\
& \cong\Ext^i_H(D^{\mu'},S^{\lambda'})
\end{align*}
where the second isomorphism follows as $D^{m(\mu')}$ is self-dual and
the third from tensoring by the sign representation, using
(\ref{Sconj}). For $l=3$ and $i=0$ see \cite[Theorem
4.6]{psqweyl}.
\end{proof}

It possible to extend the first part of the above corollary for $l=3$
to $\Ext^1$ between row three-regular three-part partitions, using
Ringel duality.  As $S(n,r)$ is quasihereditary, there is for each
$\lambda\in \partn(n,r)$ a corresponding indecomposable tilting module
$T(\lambda)$, and hence we can consider the Ringel dual of $S(n,r)$.

When $r\leq n$, Donkin has shown \cite[Proposition 4.1.4]{donkbk} that
$S(n,r)$ is Ringel self-dual (generalising the classical result first
proved in \cite{dontilt}). If $l>n$ then for all $r$ there exists
an ideal $I(n)$ in $H_{q,r}$ such that the Ringel dual of $S(n,r)$ can
be identified with $H_{q,r}/I(n)$. Further, the standard modules for
this quotient algebra are precisely the Specht modules labelled by
partitions with at most $n$ parts. (See \cite{erdqh} for the classical
case, with \cite[Section 4.7]{donkbk} for the general case.) This
induces an equivalence of categories between the category of
$S(n,r)$-modules with a filtration by induced modules and the category
of $H_{q,r}$-modules with a filtration by Specht modules labelled by
partitions with at most $n$ parts. If we restrict attention to row
regular partitions this equivalence also holds for $l=n$,
(by applying \cite[Section 4.7 (6)]{donkbk}) which
implies

\begin{propn}
If $\lambda$, $\mu \in \rpartn(3,r)$ and $l= 3$ then 
$$\Ext^1_{S(3,r)}(\nabla(\lambda), \nabla(\mu))
\cong \Ext^1_{H}(S^\lambda, S^\mu).$$
\end{propn}

It is well-known \cite{cpsk} that
\begin{equation}\label{extvan}
\Ext_{S(n,r)}^i(\nabla(\lambda), \nabla(\mu)) 
\cong
\left\{
\begin{array}{ll}
k &\mbox{if $\lambda = \mu$ and
$i=0$}\\
0 &\mbox{if $\lambda = \mu$ and
$i>0$, or if $\lambda \not < \mu$.}
\end{array}
\right.
\end{equation}

Combining this with Corollary \ref{cor:sameext} we obtain

\begin{cor}  
If $\lambda$, $\mu \in \partn(n,r)$ and $l\ge 3$ then
$$\Ext_{H}^i(S^\lambda, S^\mu) 
\cong
\left\{
\begin{array}{ll}
k &\mbox{if $\lambda = \mu$ and $i=0$}\\
0 &\mbox{if $\lambda = \mu$ and $i>0$, or if $\lambda\not <\mu$.}
\end{array}
\right.
$$
for $0 \le i \le l-3$.
\end{cor}

In particular $\Ext_{H}^1(S^\lambda, S^\lambda) =0$, if $l \ge
4$, and so this Corollary can be regarded as a generalisation of
results of Kleshchev and Nakano \cite{klenak} and Hemmer and Nakano
\cite{hemnak}.

Let $\lambda^1$ and $\mu^1 \in \partn(n_1,r_1)$ and $\lambda^2$ and
$\mu^2 \in \partn(n_2,r_2)$ be such that $\lambda^1_{n_1}\geq
\lambda^2_{1}$ and $\mu^1_{n_1}\geq \mu^2_{1}$. We can form a new
partition $\lambda\in \partn(n,r)$, where $n=n_1+n_2$ and $r=r_1+r_2$,
by defining $\lambda_i=\lambda^1_i$ if $i\leq n_1$ and
$\lambda_i=\lambda^2_{i-n_1}$ if $n_1<i\leq n_1+n_2$ (and similarly
form a new partition $\mu$). We will say that such a pair of
partitions $(\lambda,\mu)$ has a \emph{horizontal cut}.  For such a pair
we can reduce the calculation of Ext-groups between induced modules to
the same calculation for lower rank groups using

\begin{thm}[Donkin] \label{cutweyl}
Suppose that $(\lambda,\mu)$ has a horizontal cut. Then for all
 $m\geq 0$ we have
$$
\Ext^m_{S(n,r)}(\nabla(\lambda), \nabla(\mu))
\cong
\bigoplus_{m = m_1+m_2} \Bigl(
\Ext^{m_1}_{S(n_1,r_1)}(\nabla(\lambda^1),\nabla(\mu^1))
 \otimes
\Ext^{m_2}_{S(n_2,r_2)}(\nabla(\lambda^2),\nabla(\mu^2))\Bigr).
$$
\end{thm}

\begin{rem}
(i) This is a result that really goes back to \cite{donklnm} (to a
result about Levi subgroups), although it was not explicitly stated in
this form.  A more explicit version can be found in Erdmann
\cite{erd}, and in the form above in \cite{donkpreprint}.\\
(ii) We have only stated a special case of Donkin's result for type
$A$; the general version holds for any reductive group. 
\end{rem}

Combining this with Corollary \ref{cor:sameext}
we obtain

\begin{cor}[Donkin]\label{cutspecht} 
Suppose that $(\lambda,\mu)$ has a horizontal cut. Then for all $0\leq
 m\leq l-3$
 we have
\begin{align*}
\Ext^m_{H_{q,n}(k)}(S^{\lambda},S^{\mu})
\cong
\bigoplus_{m = m_1+m_2} \Bigl(
\Ext^{m_1}_{H_{q,n_1}(k)}(S^{\lambda^1},S^{\mu^1})
 \otimes
\Ext^{m_2}_{H_{q,n_2}(k)}(S^{\lambda^2},S^{\mu^2})\Bigr).
\end{align*}
\end{cor}

In the classical case a proof of this last result when $m=0$ (given
entirely in the context of the symmetric group) can be found in
\cite{fayly}. 

We say two partitions admit a \emph{vertical cut} if their conjugates admit an
horizontal cut. There are vertical cut analogues of Theorems
\ref{cutweyl} and \ref{cutspecht}, which follow from the above
results and Lemma \ref{conj}.

There is an analogue of Theorem \ref{cutweyl} in
\cite[Corollary 10]{cpstransf}.

\begin{thm}[Cline-Parshall-Scott]\label{cpscut}
Suppose that $(\lambda,\mu)$ has a horizontal cut as in Theorem
 \ref{cutweyl}. Then for all $m\geq 0$ we have
$$
\Ext^m_{S(n,r)}(L(\lambda), \nabla(\mu))
\cong
\bigoplus_{m = m_1+m_2} \Bigl(
\Ext^{m_1}_{S(n_1,r_1)}(L(\lambda^1),\nabla(\mu^1))
 \otimes
\Ext^{m_2}_{S(n_2,r_2)}(L(\lambda^2),\nabla(\mu^2))\Bigr).
$$
\end{thm}

\begin{rem} (i) As for Theorem \ref{cutweyl}, we have only stated the type
$A$ version of \cite[Corollary 10]{cpstransf}. The general version
holds for all reductive groups.\\ 
(ii) Even in type $A$, the result in
\cite{cpstransf} is actually stated rather differently, corresponding
to the form of Theorem \ref{cutweyl} given in \cite{erd}. However a
standard choice of Levi subgroups, together with an easy application
of the K\"unneth formula, gives the result in the form stated above.\\
(iii) This theorem and Theorem \ref{cutweyl} are both special cases of
a theorem of Donkin's proved in \cite[Theorem 4.2]{lylmat}.\\
(iv) Note that we cannot combine Theorem 
\ref{cpscut} and the second part of Corollary \ref{cor:sameext}
to obtain an analogue of Corollary \ref{cutspecht}, as
the Ext-groups in each case are no longer of the same form. Indeed if
we consider whether 
$
\Ext^m_{H_{q,n}}(D^\lambda, S^\mu)
$
is the same as
$$
\bigoplus_{m = m_1+m_2} \Bigl(
\Ext^{m_1}_{H_{q,n_1}}(D^{\lambda^1},S^{\mu^1})
 \otimes
\Ext^{m_2}_{H_{q,n_2}}(D^{\lambda^2},S^{\mu^2})\Bigr).
$$ 
where $(\lambda, \mu)$ has a v-cut or a h-cut,
then this is false in general, as the following counter example shows.

We know that $ \Ext^m_{H_{q,n}}(D^\lambda, S^\mu) \cong 
\Ext^m_{S(n,r)}(\nabla(\lambda'), L(\mu'))$ for $m \le l-3$, using
Corollary \ref{cor:sameext}, so we give a counter
example for algebraic groups.
Let $ l \ge 4$, $\lambda =(l-1,1)$, $\lambda^1=(2,1)$ and
$\lambda^2=(l-3)$.
So $\lambda' = (2, 1^{l-2})$, ${\lambda^{1}}' = (2,1)$ and
${\lambda^{2}}'=(1^{l-3})$. The induced module 
$\nabla{(2,1^{l-2})}$ is not simple, it has two composition 
factors ---
namely $L(1^l)$ as its head and $L(2,1^{l-2})$ as its socle.
Thus $\Hom(\nabla(\lambda'), L(\lambda')) = 0$.
But $\Hom(\nabla({\lambda^1}'), L({\lambda^1}')) = k$ and
$\Hom(\nabla({\lambda^2}'), L({\lambda^2}')) = k$ as both these induced
modules are simple and so the formula cannot hold in general.

Similarly if $\lambda'=(3, 1^{l-3})$, ${\lambda^1}'=(1^{l-2})$ and
${\lambda^2}'=(2)$ then we also get 
$\Hom(\nabla(\lambda'), L(\lambda')) = 0$,
but $\Hom(\nabla({\lambda^1}'), L({\lambda^1}')) \otimes 
\Hom(\nabla({\lambda^2}'), L({\lambda^2}')) = k$. 
\end{rem}

We present one more general type of result before we look at which
$\Ext$-groups have been explicitly calculated. This is a
generalisation of a recent result of Fayers \cite{fayers}.

We first define a function on $X^+$.  Given an integer $s$, we define
$\DD_{(s,n)}{\lambda}$ for $\lambda \in X^+$ to be the dominant weight
$(s- \lambda_n, s-\lambda_{n-1}, \ldots, s-\lambda_1)$.  (This is the
$\breve{\lambda}$ of Fayers \cite{fayers}, but we want to emphasise
the $s$ and $n$.)
Let $D=L(1,1,\ldots,1)$ be the determinant module for $\GL_n$. We
interpret $D^{\otimes s}$ to be $D$ tensored with itself $s$ times if
$s$ is positive and to be the module $L(-1,-1,\ldots,-1)\cong D^*$
tensored with itself $-s$ times if $s$ is negative. We set $D^{\otimes
0}$ be the trivial module $L(0,0,\ldots,0)\cong k$.

\begin{lem}
For all $M$, $N \in \Mod(G)$ and $i\geq 0$, we have a canonical isomorphism
$$\Ext_G^i(M,N) \cong 
\Ext_G^i((M^*)^\circ\otimes D^{\otimes s},(N^*)^\circ\otimes D^{\otimes s})
.$$
\end{lem}
\begin{proof}
We have that
$$\Ext_G^i((M^*)^\circ\otimes D^{\otimes s},
(N^*)^\circ\otimes D^{\otimes s})
\cong
\Ext_G^i((M^*)^\circ,(N^*)^\circ)$$
as the module $D^{\otimes s}$ is one-dimensional.
This latter $\Ext$ group is isomorphic to
$$\Ext_G^i(N^*,M^*)
\cong
\Ext_G^i(M,N)$$
after removing both duals.
\end{proof}

Thus we may define a functor
$\DD_{(s,n)}: \Mod(G) \to \Mod(G)$ which takes $M \in \Mod(G)$ to the module
$((M^*)^\circ)\otimes D^{\otimes s}$, and this is an equivalence on
$\Mod(G)$. The reason for abusing notation and reusing the notation
$\DD_{(s,n)}$ is made clear by the next lemma.

\begin{lem}\label{Ddoes}
Let $s \in \Z$, and $\lambda$, $\mu \in X^+$. Then
\begin{enumerate}
\item[(i)]{
$\nabla(\DD_{(s,n)}(\lambda)) \cong \DD_{(s,n)}(\nabla(\lambda))$
}
\item[(ii)]{
$\Delta(\DD_{(s,n)}(\lambda)) \cong \DD_{(s,n)}(\Delta(\lambda))$
}
\item[(iii)]{
$L(\DD_{(s,n)}(\lambda)) \cong \DD_{(s,n)}(L(\lambda))$
}
\item[(iv)]{
$T(\DD_{(s,n)}(\lambda)) \cong \DD_{(s,n)}(T(\lambda))$
}
\end{enumerate}
\end{lem}
\begin{proof}
We prove (i), the other statements are similar.
We have $\nabla(\lambda) ^* \cong \Delta(-w_0\lambda)$ where $w_0$ is
the longest element of the associated Weyl group $\Sym_{n}$. If $\lambda
= (\lambda_1, \lambda_2, \ldots, \lambda_n)$ then
$-w_0\lambda = (- \lambda_n, -\lambda_{n-1}, \ldots,-\lambda_1)$.
We thus have $(\nabla(\lambda)^*)^\circ \cong
\nabla(\DD_{(0,n)}(\lambda)).$
Hence $\DD(\nabla(\lambda)) \cong \nabla(\DD_{(0,n)}(\lambda)) \otimes
D^{\otimes s} \cong \nabla(\DD_{(s,n)}(\lambda))$.
\end{proof}

\begin{cor}\label{Dgives}
For all $\lambda$, $\mu\in X^+$, $s\in\Z$, and  $i\geq 0$ we have
$$\Ext_{G}^i(A(\lambda),B(\mu)) \cong 
\Ext_{G}^i(A(\DD_{(s,n)}(\lambda)),
B(\DD_{(s,n)}(\mu)))$$
where $A(-)$ and $B(-)$ are any of $\nabla(-)$, $\Delta(-)$, $L(-)$
and $T(-)$. 
\end{cor}

When all the weights in the above corollary are polynomial we get the
same isomorphism for the $q$-Schur algebra.  Using Corollary
\ref{cor:sameext} we obtain the following result for the corresponding
Hecke algebra.

\begin{thm}
Suppose that $\lambda$, $\mu \in \partn(n,r)$ and that $s$ is greater
than $\lambda_1$ and $\mu_1$.  Then for $l \ge 3$ and $0 \le i \le
l-3$ we have that
$$\Ext_H^i(S^\lambda, S^\mu) \cong \Ext_H^i(S^{\DD_{(s,n)}(\lambda)},
S^{\DD_{(s,n)}(\mu)})$$ 
and
$$\Ext_H^i(D^{\lambda'}, S^{\mu'}) \cong \Ext_H^i(D^{\DD_{(s,n)}(\lambda')},
S^{\DD_{(s,n)}(\mu')})$$ 
We also get the first isomorphism for $i=0$
and $l=2$ when either $\lambda$ is row $2$-regular or $\mu$ is
column $2$-regular.
\end{thm}
\begin{proof}
We combine Theorem \ref{thm:donk} and Theorem \ref{thm:equalhom} 
and apply them to the previous corollary.
(Note that $\lambda$ is row $l$-regular if and only if
$\DD_{(s,n)}(\lambda)$ is row $l$-regular)
\end{proof}

\begin{rem}
(i) In the case $i=0$ and $p \ge3$, a proof of the first part of this
result located entirely in the symmetric group setting has been given
by Fayers \cite{fayers}. It is also known to be false in full generality
for $p=2$. \\ 
(ii) We have only stated the type $A$ version of the
equivalence of categories given by $\DD_{(s,n)}$, and of Corollary
\ref{Dgives}. The general version holds for all reductive groups.\\
(iii) In type $A$, this equivalence has been realised via an
isomorphism of generalised Schur algebras by Feng, Henke, and K\"onig
\cite{fhk}. 
\end{rem}

\section{Calculating Ext-groups I: Homomorphisms}\label{homsec}

In the remaining sections we will consider various results which give
explicit values for certain $\Ext$-groups. We begin with the
homomorphism case.  Very little is known about homomorphisms between
Weyl (or Specht) modules in general. Almost all the general results
are either for weights which are \lq close together\rq\ for example
those of \cite{kop1}, or those where the calculation can be reduced to
a low rank case using Theorem \ref{cutweyl}. The only exception to
this, due to Carter and Payne, still only considers weights which are
related by a single reflection (in the alcove-theoretic sense). As
many of the results in this section are most naturally expressed in
this geometric language we will switch freely between the partition
and alcove contexts.  In this section we take $q=1$ unless
explicitly indicated otherwise.

We say that a weight $\mu$ is a \emph{Steinberg} weight if $\mu_i - \mu_{i+1}
\equiv -1 \pmod p$ for all $i$. We let $W_p$ be the affine Weyl group
for $G$ which acts on $X$ via the dot action. This defines a system of
hyperplanes and facets for $X$: details may be found in
\cite[II, Section 1.5 and Chapter 6]{jantzen}. Two dominant weights
can only be in the same block if they are in the same $W_p$-orbit.
Koppinen \cite[Theorem 7.1]{kop1} proves the following.
\begin{thm}[Koppinen]\label{kopthm}
Let $\mu$ be a dominant Steinberg weight, and $W^{\mu}_p$ be the
stabiliser of $\mu$ in $W_p$. Take $\lambda$ lying in a facet whose
closure contains $\mu$ and set $J = W^{\mu}_p \cdot \lambda$.
Then\\ 
\noindent(i)
If $\xi$, $\xi' \in J$ and $\xi \ge \xi'$ then
$\Hom_G(\nabla(\xi), \nabla(\xi'))\cong k$.
\\
\noindent(ii)
If $\xi$, $\xi'$, $\xi'' \in J$ and $\nabla(\xi) \to \nabla(\xi') \to
\nabla(\xi''')$ are non-zero homomorphisms then the composite map is
nonzero.
\end{thm}

This theorem has been slightly generalised by Wen \cite{wen}, who
also has more to say about the higher $\Ext$ groups. We will return to
these results in section \ref{higher}.

Let $\lambda, \mu \in \partn$. Suppose that there exists a pair  $i >
j$, and  integers $n>0$ and $m$ such that $\mu$ is obtained from
$\lambda$ by moving $d$ boxes from the $i$th row to the $j$th row,
with $d = \lambda_i -\lambda_j -i +j -m p^a < p^a$. In this case we
say that $\mu$ is a \emph{local $p^a$-reflection of $\lambda$}.

We have the following theorem due to Carter and Payne \cite{cpmaps}.

\begin{thm}[Carter-Payne]\label{cpmap}
Let $\mu$ be a local $p^a$-reflection of $\lambda$ for some $a\geq 0$.
Then we have
$$\Hom_{S(n,r)}(\nabla(\lambda), \nabla(\mu)) \ne 0$$
and if $p>2$ then 
$$\Hom_{k\Sym_r}(S^\lambda, S^\mu) \ne 0.$$
\end{thm}

This is proved first for $\SL_n$ by looking at the action of the
hyperalgebra on Weyl modules, and then by applying Theorem
\ref{thm:equalhom} to get the result for the symmetric group. It has
recently been generalised by Fayers and Martin \cite{faymar}, who work
in the symmetric group setting by looking at semi-standard
homomorphisms. We will return to this generalisation in Section \ref{FMsec}.

\begin{rem}
(i) If $j = i+1$ then the corresponding quantum result also holds, but
  replacing occurences of $p^a$ by $lp^{a-1}$. Further, if $j=i+1$
  then these are the only homomorphisms for the Schur algebra, and for
  the Hecke algebra if $l \ge 3$. Thus in both cases the $\Hom$-spaces
  are at most one-dimensional. This result is well-known; a proof can
  be found in \cite[Theorem 5.1]{coxerd}.\\ 
(ii) The reflection terminology is motivated by the alcove-geometric
  approach to weights. Using the standard dot action of the Weyl group
  for $G$, this condition corresponds to $\mu$ being the reflection of
  $\lambda$ about some $p^a$-hyperplane $P$ such that no parallel
  $p^a$-hyperplane lies between $\lambda$ and $P$. \\
(iii) The condition $p>2$ is necessary for the symmetric group result.
  For example when $p=2$ we have that $S^{(1,1)}$ is isomorphic to
  $S^{(2)}$ and so the corresponding Hom-space will be non-zero.\\
(iv) We expect the corresponding quantum result to hold. However the
  proof in the classical case requires explicit manipulation of the
  action of the hyperalgebra, which will be much more complicated in the
  quantum case. The quantum result is known to hold if $j=i+1$ as in
  (i) above or if $j=i+2$, see \cite[theorem 9.1]{par4} and apply
  Theorem \ref{cutweyl}.
\end{rem}

Kulkarni \cite{kul1, kul2} has calculated  certain
Ext groups between Weyl modules over the integers, and his results
have consequences for the modular theory. In particular, he shows that
$$\dim\Hom(\nabla(r,0,\ldots,0), \nabla(\lambda))\leq 1$$
and can describe precisely when the Hom-space is non-zero
\cite[Theorem 2.2 and the following Remark (3)]{kul1}.

Apart from the $\GL_2$-result mentioned above, the only other
classical case
where homomorphisms have been completely determined is for $\GL_3$,
(with $p>2$) where we have recently given a complete classification
\cite{coxpar}. This is given by a complicated recursive set of
interlocking results, expressed in the language of alcove geometry,
which is too lengthy to reproduce here. Instead we will outline the
main features of the classification.

The basic idea is to use certain nice filtrations of induced modules,
called good $p$-filtrations. These exist for all induced modules, and
for $\GL_3$ the structure of these filtrations has been completely
determined in \cite{par1}. A \emph{good $p$-filtration} is one where
successive quotients are all of the form
$\nabla_p(\mu)=\nabla(\mu'')^F\otimes L(\mu')$ where
$\mu=\mu'+p\mu''$ for some dominant weight $\mu''$ and $p$-restricted
weight $\mu'$. Given a pair of such modules $\nabla_p(\mu)$ and
$\nabla_p(\tau)$  we have 
\begin{equation}
\label{twist}
\Hom_G(\nabla_p(\mu),\nabla_p(\tau))\cong\left\{\begin{array}{ll}
\Hom_G(\nabla(\mu''),\nabla(\tau''))\quad &\text{if $\mu'=\tau'$}\\
0 &\text{otherwise}\end{array}\right.
\end{equation}
which enables us to proceed by induction.

\begin{figure}[ht]
\centerline{\epsffile{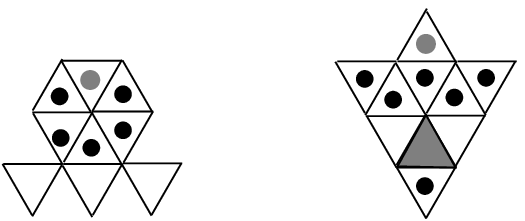}}
\caption{\label{alcoves}}
\end{figure} 

When the weight $\lambda$ is small then all terms in the good
$p$-filtration of $\nabla(\lambda)$ are simple, and the homomorphisms
from $\nabla(\lambda)$ can be classified using the results in
\cite{par1}. For example, suppose that $\lambda$ lies in an alcove
labelled by a lightly shaded dot in Figure \ref{alcoves}, and suppose
it is such that all nine unshaded alcoves in the respective diagram
are dominant. Then these alcoves label the terms in the good
$p$-filtration for $\nabla(\lambda)$, and the marked alcoves indicate
the possible homomorphisms in the case where these terms are all simple.

In the general case we use translation functor arguments to translate
from a Steinberg weight $\theta$ just below $\lambda$ --- where we
know all homomorphisms by induction, as
$\nabla(\theta)\cong\nabla_p(\theta)$. (If no such $\lambda$ exists we
use an alternative argument based on an identification of certain
induced modules with symmetric powers.) Each homomorphism from
$\nabla(\theta)$ to some $\nabla(\tau)$ gives rise to the possibility
of a map from $\nabla(\lambda)$ into certain induced modules labelled
by weights \lq near' to $\tau$. This, together with the known block
structure of the module category, provides a necessary condition for
there to be a non-zero homomorphism from $\nabla(\lambda)$ to
$\nabla(\mu)$. As $L(\mu)$ is the socle of $\nabla(\mu)$ another
obvious necessary condition is that $L(\mu)$ is a composition factor
of $\nabla(\lambda)$.

It is almost always the case that this pair of necessary conditions is
also sufficient. We can also give an explicit description of the
exceptions to this rule. All these cases are most easily described
graphically in a series of diagrams similar in nature to those in
Figure \ref{alcoves}, with one such diagram for each homomorphism from
$\nabla(\theta)$.

To show that these conditions are sufficient, it is enough to
construct homomorphisms in each of the remaining cases. There are two
obvious classes of homomorphisms already at our disposal.  The first
of these are the Carter-Payne maps coming from Theorem \ref{cpmap},
and we begin by determining a necessary condition for composites of
such maps to be non-zero. The second obvious class consists of those
from the top term in the good $p$-filtration (which is uniquely
defined). By (\ref{twist}) these are twists of maps from induced
modules labelled by smaller weights, which are assumed known by
induction.

By considering all the various possible cases we can show that almost
all non-zero homomorphisms between induced modules are either
composites of Carter-Payne maps or are twisted maps. Unfortunately,
there are again exceptions to this rule, arising from certain explicit
configurations of weights which we can completely classify. These
exceptional maps can be constructed from pairs of maps on terms in
the good $p$-filtration.

At several key stages during the proof we need the following 
result, which is also proved as part of the induction argument

\begin{thm} Suppose that $p>2$. For any pair $\lambda$, $\mu$ of
  dominant weights for $\GL_3$, we have
$$\dim\Hom(\nabla(\lambda),\nabla(\mu))\leq 1.$$
\end{thm}

Putting this together with the argument outlined above, and Theorem
\ref{thm:equalhom} we obtain

\begin{thm} Suppose that $p>2$. For any pair $\lambda$, $\mu$ of
  dominant weights for $\GL_3$, we can classify all homomorphisms between
$\nabla(\lambda)$ and $\nabla(\mu)$ . If $\lambda$ and $\mu$ are both
  polynomial weights then we can also classify all homomorphisms between
$S^\lambda$ and $S^\mu$.
\end{thm}

\begin{rem}
(i) A similar result holds for the quantum case. The main obstacle to
  proving this in \cite{coxpar} was the absence of quantum versions of
  the good $p$-filtration structure theorems from \cite{par1} and of
  the Carter-Payne theorem. The former, and a special case of the
  latter sufficient for this problem, have now been proved by the
  second author \cite{par4}. \\ (ii) It is hard to see how the result
  for $\GL_3$ could be stated using combinatorics of partitions. We
  believe that the alcove-theoretic formalism has been under-utilised
  in the Hecke algebra setting; further applications of this formalism
  can be found in \cite{cox5}.\\ (iii) One might hope that all
  Hom-spaces were at most one dimensional, by analogy with the Verma
  module case. However, while we know of no examples of Hom-spaces for
  induced modules which are greater than one dimensional, this may
  simply be because such Hom-spaces are only known in certain very
  special cases. The methods used to date in such calculations provide
  no compelling indication either way.  As we will see in Remark
  \ref{fmrem}(iii), a potential source of large Hom-spaces has been
  given by Fayers and Martin \cite{faymar}; however no actual examples
  have yet been found.
\end{rem}

\section{A generalisation of the Carter-Payne theorem}\label{FMsec}

In this section we assume that $q=1$.  Fayers and Martin \cite{faymar}
have given a sufficient condition for the existence of a non-zero
homomorphism from $S^{\lambda}$ to $S^{\mu}$ when $\mu$ is obtained
from $\lambda$ by moving $s$ boxes from one row of $\lambda$ to a
lower row. The condition is rather a complicated one; we will review
the basic combinatorial framework needed to state the result, and then
indicate a few simplifications that can be made. (In a similar vein,
K\"unzer \cite{kun1,kun2} has considered the cases $s=1$ and $2$,
giving an explicit construction of a morphism.) In what follows we
assume standard facts about tableaux etc.

Let $T$ be a $\lambda$-tableau of weight $\mu$ (i.e. shape $\lambda$
and with $\mu_i$ entries equal to $i$). For a tableau $T$ let $T_j^i$
be the number of entries equal to $i$ in row $j$, and
$T_h=\sum_{j>h}T_h^j$. Given integers $a$ and $b$ we set $a^{b
\downarrow}=\prod_{i=0}^{b-1}(a-i)$ and $a^{b
\uparrow}=\prod_{i=0}^{b-1}(a+i)$.

Using Corollary \ref{cutspecht} , we may reduce to the case where
$$\lambda=(l_0+s,l_1^{m_1-1},l_2^{m_2-m_1},\ldots,l_r^{m_r-m_{r-1}})
\quad\wand\quad
\mu=(l_0,l_1^{m_1-1},l_2^{m_2-m_1},\ldots,l_r^{m_r-m_{r-1}},s).$$ Pick
$\gamma_i\in\Z$ for $1\leq i\leq r$, and set $m_0=1$. With these
assumptions we define $n_i(T)=\sum_{m_{i-1}+1}^{m_i}T_h$ and

$$c_i(T)=\gamma_i^{(s-n_i(T))\downarrow}\prod_{h=m_{i-1}+1}^{m_i}T_h!.$$

We will need two functions $f$ and $g$. The first is easy to define: 
$$f(T)=\prod_{i=1}^rc_i(T).$$
Given a partition $\mu$, and integers $d\geq 1$ and $0\leq t\leq
\mu_{d+1}$ we define a new {\it composition} $\nu=\nu(\mu,d,t)$ by
setting 
$$\nu_i=\left\{\begin{array}{ll} \mu_i+t&\wif i=d\\ 
\mu_i-t&\wif i=d+1\\ 
\mu_i&\otherwise.\end{array}\right.$$ 
Given $\lambda$, we define $g_{\mu,b,t}(S)$ for a $\lambda$-tableau $S$ of
type $\nu(\mu,m_b,t)$ by
\begin{equation}
\label{gis}
g_{\mu,b,t}(S)=\left\{\begin{array}{ll}
f(S)\frac{(l_0-l_1+m_1-1+s-\gamma_1)^{t\uparrow}}{t!}& \wif
b=0\\ 
\frac{f(S)}{c_b(S)}
\frac{(l_b-l_{b+1}+m_{b+1}-m_{b}+\gamma_b-\gamma_{b+1})^{t\uparrow}}{t!}
(\gamma_b^{(s-n_b(S)-t)\downarrow}(S_{m_b}+t)!)\prod_{j=m_{b-1}+1}^{m_b-1}S_j!&
\wif 0<b<r\\
\frac{f(S)}{c_b(S)}
\frac{(l_r-s+1+\gamma_r))^{t\uparrow}}{t!}
(\gamma_b^{(s-n_b(S)-t)\downarrow}(S_{m_b}+t)!)\prod_{j=m_{b-1}+1}^{m_b-1}S_j!&
\wif b=r.\end{array}\right.
\end{equation}

Finally, we say a tableau is \emph{pseudo-standard} if the entries
weakly increase along rows and $T_j^i$ is non-zero only if $i=j$ or
$\lambda_i<\lambda_j$. A \emph{nice} tableau is one where the number of
entries greater than $m_i$ in the set of all rows of length $l_i$ is
at most $l_i-l_{i+1}$, for each $i$.

The following result is due to Fayers and Martin \cite[Theorem
  22]{faymar}, and (although this is not obvious) generalises Theorem
  \ref{cpmap}.

\begin{thm}[Fayers-Martin]\label{mf} Suppose that $\lambda$ and $\mu$ are as above, and that
 we can pick integers $e>0$ and $\gamma_i$ such that:
\begin{itemize}
\item[(i)] For some nice pseudo-standard $\lambda$-tableau $T$
  of type $\mu$ the coefficient $f(T)$ is not divisible by $p^e$.
\item[(ii)] For all $0\leq i\leq r$ and $1\leq t\leq \mu_{m_i+1}$, and all
  pseudo-standard $\lambda$-tableaux $S$ of type 
  $\nu(\mu,m_i,t)$, the coefficient $g_{\mu,i,t}(S)$ is divisible by $p^e$.
\end{itemize}
Then there
  exists a non-zero homomorphism from $S^{\lambda}$ to $S^{\mu}$.
\end{thm}

It will be useful to simplify the expression for $g_{\mu,b,t}(S)$.
Let $\gamma_0=s$, $m_{r+1}=m_r+1$, $l_{r+1}=s$, and
$\gamma_{r+1}=0$. Then the expression above simplifies to
\begin{equation}
\label{galt}
g_{\mu,m_b,t}(S)=f(S)
\frac{(l_b-l_{b+1}+m_{b+1}-m_{b}+\gamma_b-\gamma_{b+1})^{t\uparrow}}{t!}
(S_d+t)^{t\downarrow}
\frac{1}{(\gamma_b+n_b(S)+t-s)^{t\downarrow}} 
\end{equation}
for all $0\leq b\leq r$.

This theorem is rather difficult to work with. As has already been
noted, it gives a generalisation of Theorem \ref{cpmap}, but even this
is not clear. (Fayers and Martin do not give Theorem \ref{cpmap} as a
direct consequence of their main result, but instead show that it
can be proved by methods which are then generalised to give Theorem \ref{mf}.)

There are many freely chosen variables, and for a given
$\lambda$ and $\mu$ potentially many $g$'s to calculate. It is not
even clear as stated where the requirement that $\lambda$ and $\mu$
are in the same block will play a role. We will show that the choice
of the $\gamma_i$ is essentially unique (modulo $p$), and how the
block restriction then appears in the combinatorics.

\begin{propn}\label{notfree}
 Suppose that for given $\lambda$ and $\mu$ the integers
  $\gamma_i$ are such that conditions (i) and (ii) of Theorem \ref{mf}
  hold. Then we must have
$$\gamma_i\equiv l_0-l_i+m_i+s-1 \pmod p$$
for $1\leq i\leq r$, and 
$$l_0+m_r\equiv 0 \pmod p.$$
\end{propn}
\begin{proof}
 First suppose that $d=1(=m_0)$ and that $T$ is as in
Theorem \ref{mf}. We will take $t=1$ and construct a new tableau
$S$. If there is a 2 in the first row of $T$ let $S$ be the tableau
obtained by replacing this with a $1$; then $f(S)=f(T)$. By assumption
we know that $p^e$ divides $f(S)$ but not $f(T)$, and hence from
(\ref{gis}) we deduce that 
$$\gamma_1\equiv l_0-l_1+m_1-1+s \pmod p$$
as required.

If there is not a $2$ in the first row of $T$ then
there must be some entry greater than $2$ in the first row: exchange
this for a $2$ from the second row (which must contain only $2$s) and
then replace the $2$ by a $1$. In this case we have $n_1(S)=n_1(T)+1$
and $S_2!=1!=0!=T_2$. Therefore 
$$f(T)=(\gamma_1-s+n_1(T))f(S)$$ 
and from (\ref{gis}) we deduce that
$$g(S)(\gamma_1-s+n_1(T))=f(T)(l_0-l_1+m_1-1+s-\gamma_1).$$
As $p^e$ divides $g(S)$ but not $f(T)$ we must have that $p$ divides 
$(l_0-l_1+m_1-1+s-\gamma_1)$, as required.

The cases $d>1$ are similar; we set $t=1$ and indicate how to
construct a suitable $S$ from the initial choice of $T$.

Let $d=m_b$ with $1\leq b<r$. We will proceed by induction on $b$; the
case $b=0$ having already been considered, and hence assume that
$\gamma_b$ has already been determined. First suppose that
$T_d>0$. Then there exists a $d$ in $T$ in a row above row $d$ and an
entry $a>d$ at the end of row $d$. If there exists a $d+1$ in one of
the first $d$ rows then replace it by $d$ and swap this new entry with
$a$ (if it is not already the entry $a$)to form $S$. In this case
$S_d=T_d-1$ and $n_b(S)=n_b(T)-1$ and we have
\begin{equation}\label{midcase}
f(S)=f(T)\frac{(\gamma_b-s+n_b(T))}{T_d\delta}
\end{equation}
where $\delta=1$. If there is not a $d+1$ in one of the first $d$ rows
then $T_{d+1}=0$, and we let $S$ be the tableau obtained from $T$ by
replacing the final element of row $d+1$ by $d$ and swapping it with
the $a>d+1$ at the end of row $d$. Then $S_d=T_d-1$ with
$n_b(S)=n_b(T)-1$, and $S_{d+1}!=T_{d+1}!=1$ with
$n_{b+1}(S)=n_{b+1}(T)+1$ and $f(S)$ is given by equation
(\ref{midcase}) where $\delta=\gamma_{b+1}-s+1+n_{b+1}(T)$.

Thus when $T_d>0$ we have from (\ref{midcase}) and (\ref{galt}) that 
$$g(S)=
f(T)
\frac{(l_b-l_{b+1}+m_{b+1}-m_{b}+\gamma_b-\gamma_{b+1})
(S_d+1)(\gamma_b+n_b(S)+1-s)}{(S_d+1)\delta(\gamma_b+n_b(S)+1-s)}
$$
and hence 
$$\delta g(S)=f(T)(l_b-l_{b+1}+m_{b+1}-m_{b}+\gamma_b-\gamma_{b+1}).$$ 
Arguing as in the $d=1$ case  we deduce that 
$$l_b-l_{b+1}+m_{b+1}-m_{b}+\gamma_b-\gamma_{b+1}\equiv 0  \pmod p$$
which by induction implies that $\gamma_{b+1}$ is of the required form.

If $T_d=0$ then there exists some entry $a\geq d+1$ in a row above row
$d$. If $T_{d+1}>0$ then we may chose $a=d+1$, and $S$ is obtained by
replacing $a$ by $d$, and $f(S)=f(T)$ (and $S_d+1=1$). If $T_{d+1}=0$
then change the final $d+1$ to a $d$ and swap this entry with the
entry $a$. In either case it is now easy to verify as above that
$\gamma_{b+1}$ must be of the required form.

Finally, suppose that $d=m_r$ (and so $\lambda$ has no row $d+1$). If
there exists a $d+1$ in a row above row $d$ in $T$ let $S$ be obtained
by replacing this with a $d$. Otherwise there must be a $d+1$ in row
$d$; again replace this by a $d$. In either case by arguing as above
we see that we must have 
$$\gamma_r+1\equiv l_0-l_{r+1}+m_{r+1}+s-1 \pmod p.$$
But we fixed $l_{r+1}=s$, $\gamma_{r+1}=0$ and $m_{r+1}=m_r+1$, which
implies that 
$$l_0+m_r\equiv 0 \pmod p$$ 
as required.
\end{proof}

\begin{rem}(i) Suppose that $\lambda$ and $\mu$ are both partitions of
  $a$ with at most $n$ parts (for some $a$ and $n$), such that the
  conditions in Theorem \ref{mf} are satisfied.  Recall the
description of the blocks of the symmetric group: $D^\lambda$ and
$D^\mu$ are in the same block of $k\Sigma_a$ if and only if there
exists a permutation $\sigma\in\Sigma_n$ such that 
$$\lambda_i-i\equiv \mu_{\sigma(i)}-\sigma(i) \pmod p.$$ By our
initial assumptions on $\lambda$ and $\mu$ we have $\lambda_i=\mu_i$
for $1<i<m_r$ and hence to be in the same block we must have either
$\lambda_i\equiv \mu_i \pmod p$ for all $i$ or $\lambda_1-1\equiv
\mu_{m_r+1}-(m_r+1) \pmod p$. But this latter condition becomes
$l_0+s-1\equiv s-m_m-1 \pmod p$, which is automatically satisfied by
Proposition \ref{notfree}.\\ (ii) Fayers and Martin have already
observed \cite[Remark after Theorem 22]{faymar} that their result is
most useful when a stronger version of the conditions in this
Proposition hold: taking all but one of the congruences for the
$\gamma_i$ to be equalities. This observation was based on empirical
evidence; our result provides theoretical corroboration for this.\\
\end{rem}

Unfortunately, the conditions in Theorem \ref{mf} are very difficult
to work with, even after introducing the simplifications in
Proposition \ref{notfree}. We conclude this section with a brief
discussion of what is known, and what might be hoped to be true.

\begin{rem}\label{fmrem}
(i) The only new example given in \cite{faymar} as an application of
their main theorem is the case $\lambda=(7,3)$ and $\mu=(4,3^2)$, with
$p=3$. This case would be covered if we could relax the condition on
$d$ in the definition of local $p^a$-reflection from $d<p^a$ to $d\leq
p^a$. (Of course there would have to be further conditions present, as
we cannot relax the definition in this way for two-part partitions.)
It would be interesting to know whether there are any examples given
by Theorem \ref{mf} which are not of this form.  \\ (ii) There
certainly are more general pairs of weights related by a single
reflection than would be obtained simply by relaxing the above
inequality; an examples for three-part partitions is given by the 
homomorphism from the highest to lowest alcove in the righthand diagram
in Figure \ref{alcoves}. As stated Theorem \ref{mf} gives only a
sufficient condition for the existence of a homomorphism; it would be
interesting to know whether it is also necessary.\\ (iii) Fayers and
Martin have also provided a potential source of greater than
one-dimensional Hom-spaces. They introduce the notion of a {\it good}
quasi-standard tableaux, and prove that homomorphisms corresponding to
such tableaux from Theorem \ref{mf} are linearly independent
\cite[Proposition 23]{faymar}. However, we know of no examples of
pairs of such tableaux --- and of course they cannot exist for any of
the examples where Hom-spaces have already been calculated.
\end{rem}

\section{Calculating Ext-groups II: The $q$-$\GL_2$ case}

In this section we will consider the problem of determining all Ext
groups between pairs of induced modules for $q$-$\GL_2$. This is the
only case in which a complete answer can be given. 

The first work in this area was by Erdmann \cite{erd}, who determined
all the $\Ext^1$ between induced modules for $\GL_2$. 
This result was quantised by the first
author \cite{cox}, who with Erdmann refined the methods used in order
to determine $\Ext^2$  between induced modules \cite{coxerd}. 
This work was based on
computations for the corresponding Frobenius kernel (where everything
can be determined relatively easily), followed by an elementary
application of the Lyndon-Hochschild-Serre spectral sequence to transfer
these results to the algebraic group. Similar methods were used by De
Visscher \cite{deV1} to calculate $\Ext^1$ between twisted tensor
products of Weyl modules and induced modules.

Recent work of the second author \cite{par3} 
has determined $\Ext^i$ for all $i$ between induced modules
for $q$-$\GL_2$ (including the classical case $q=1$), as well as $\Ext^i(\nabla(\lambda), L(\mu))$,
$\Ext^i(L(\lambda), \nabla(\mu))$
and $\Ext^i(L(\lambda), L(\mu))$, for $\lambda$, $\mu \in X^+$.
By Corollary \ref{cor:sameext} this result for induced modules also calculates
the $\Ext^i$ between Specht modules for 2-part partitions provided
that $i\leq l-3$. We will describe the main features of this result.

The main idea is to refine the method first used by Erdmann, and
consider in detail the spectral sequence involved. By giving an
alternative derivation of this sequence, it becomes possible to
describe in detail various pages in the sequence. From this can be
deduced various reduction theorems for Ext-groups, which allow any
Ext-group to be computed by induction on the weights.

In order to state the main results explicitly we focus for the moment
on two part partitions.
We will distinguish the modules for classical $\GL_2(k)$ and the
quantum group $q$-$\GL_2(k)$ by putting a bar on the modules for the 
classical groups.  The two are related by the Frobenius morphism
$\mathrm{F}: q$-$\GL_2(k) \to \GL_2(k)$, and
we have an map of module categories which takes a module for
classical $\GL_2(k)$ to a module for the quantum group $q$-$\GL_2(k)$,
namely the map that sends $\overline{M}$ to the twisted module $\overline{M}\frob$.

\begin{thm}
For $a \ge b$ with $a-b$ odd, $0 \le i \le l-2$ and $m \in \N$ we
have
\begin{align*}
\Ext&_{q-\GL_2(k)}^m(\nabla(la+i,0),\nabla(lb+l-2-i+d,d))\\
&\cong \Ext_{q-\GL_2(k)}^{m-1}(\nabla(la-1,i+1), \nabla(lb+l-i-2+d,d))
\oplus
\Ext_{\GL_2(k)}^m(\bnabla(a-1,0), \bnabla(b+f,f))
\end{align*}
where $\Ext^{-1}$ is interpreted as the zero module,
$f=\frac{a-b-1}{2}$ and $d = lf +i+1$.
\end{thm}
\vskip10pt

\begin{thm}
For $a \ge b$ with $a-b$ even, $0 \le i \le l-2$ and $m \in \N$ we
have
$$
\Ext_{q-\GL_2(k)}^{m}(\nabla(la+i,0), \nabla(lb +i+d,d))
\cong \Ext_{q-\GL_2(k)}^{m}(\nabla(l(a-b)+i,0), \nabla(i+d,d))
$$
where $d = l(\frac{a-b}{2})$.
If $m\ge 1$ then also 
$$
\Ext_{q-\GL_2(k)}^{m}(\nabla(l(a-b)+i,0), \nabla(i+d,d))
\cong \Ext_{q-\GL_2(k)}^{m-1}(\nabla(l(a-b)-1,i+1), \nabla(i+d,d)).
$$
\end{thm}

\vskip10pt
\begin{thm}
For $a \ge b$ with $a-b$ even  and $m \in \N$ we
have
$$
\Ext_{q-\GL_2(k)}^m(\nabla(la+l-1,0),\nabla(lb+l-1+d,d))\\
\cong
\Ext_{\GL_2(k)}^m(\bnabla(a,0), \bnabla(b+f,f))
$$
where $f=\frac{a-b}{2}$ and $d = lf$.
\end{thm}

For simplicity we have stated these results in terms of a weight of the
form $(c,0)$. We have for all $b>0$ that 
$$\nabla(a,b)\cong \nabla(a-b,0)\otimes{\det}^b$$ 
(where $\det$ is the one dimensional module corresponding to the
determinant), and so all other cases can be
reduced to this one.  By block considerations and equation
\eqref{extvan} any non-zero Ext-group can be reduced to one of the
forms considered in the above Theorems, and hence (as we already know
the $m=0$ results and using Theorem \ref{ryom}) any Ext-group can be
computed by induction. For $m\leq 2$ explicit closed forms for these
results were given in \cite{coxerd}. Even for such small degrees the
closed forms quickly become unmanageable, and in general a recursive
formula seems to be the best that can be expected.

A recent calculation of Erdmann and Hannabuss \cite{erdhan} using the
result for $p=2$ has shown that in this case the dimension of the
$\Ext$ group $\Ext^m_{\GL_2(k)}(\nabla(2a,0), \nabla(a,a))$ is the 
same as the number of
partitions $(\eta_1, \eta_2, \ldots, \eta_{m+1})$ with $2^{\eta_1}
+2^{\eta_2} + \cdots +2^{\eta_{m+1}} - 1 = a$.

We also have
\begin{thm}
For $a \ge b$ with $a-b$ odd, $0 \le i \le l-2$ and $m \in \N$ we
have
\begin{align*}
\Ext&_{q-\GL_2(k)}^m(\nabla(la+i,0),L(lb+l-2-i+d,d))\\
&\cong \Ext_{q-\GL_2(k)}^{m-1}(\nabla(la-1,i+1), L(lb+l-i-2+d,d))
\oplus
\Ext_{\GL_2(k)}^m(\bnabla(a-1,0), \bar{L}(b+f,f))
\end{align*}
where $\Ext^{-1}$ is interpreted as the zero module,
$f=\frac{a-b-1}{2}$ and $d = lf +i+1$.
\end{thm}
\vskip10pt

\begin{thm}
For $a \ge b$ with $a-b$ even, $0 \le i \le l-2$ and $m \in \N$ we
have
$$
\Ext_{q-\GL_2(k)}^{m}(\nabla(la+i,0), L(lb +i+d,d))
\cong \Ext_{q-\GL_2(k)}^{m-1}(\nabla(la-1,i+1), L(lb+i+d,d))
$$
where $d = l(\frac{a-b}{2})$.
\end{thm}
and
\begin{thm}
For $a \ge b$ with $a-b$ even, and $m \in \N$ we
have
$$
\Ext_{q-\GL_2(k)}^{m}(\nabla(la+l-1,0), L(lb +l-1+d,d))
\cong \Ext_{q-\GL_2(k)}^{m}(\bnabla(a,0), \bar{L}(b+f,f))
$$
where 
$f=\frac{a-b}{2}$ and $d = lf$.
\end{thm}

\section{\label{higher}Calculating Ext-groups III: Higher ranks}

In this section we will complete our survey of general results for
Ext-groups between pairs of induced or Specht modules. The results
here can all be regarded as corresponding to cases where the pair of
indexing weights $\lambda$ and $\mu$ are \lq close together\rq.  To
make this notion more precise we introduce the following function of
$\lambda$.

Given $\lambda \in \partn(n,r)$ we define
$$d(\lambda)= \sum_{i=1}^{n-1}\sum_{j=i+1}^{n}
               \left\lfloor \frac{\lambda_i -\lambda_j
                        -i+j-1}{l}\right\rfloor.$$

We will consider two weights $\lambda$ and $\mu$ to be close if the
difference $|d(\lambda)-d(\mu)|$ is small.

In what follows, most of the results were originally stated using
alcove combinatorics and the standard length function on elements of
the corresponding affine Weyl group. The correspondence between this
and the partition combinatorics together with the function $d$ is
given in \cite[Lemmas 3.5 and 5.1]{par2}.

The first part of the following result is due to Andersen
\cite[Proposition 3.5]{andGBcoh} and the second to Ryom-Hansen
\cite[Theorem 2.4]{rhcrelle}. (Both stated and proved in the classical
case, but the quantum version follows in the same way.)
In each case the cited results are for
the case where one module is simple; the statements involving two
induced modules follows easily from this, confer \cite[Proposition
4.6]{par2}.

\begin{thm}[Andersen, Ryom-Hansen]\label{ryom}
Suppose $\lambda$, $\mu \in \partn$.
If $m >d(\lambda) - d(\mu)$ we have 
$$
\Ext^{m}_{S(n,r)}
(L(\lambda), \nabla(\mu)) \cong 
\Ext^{m}_{S(n,r)}
(\nabla(\lambda), \nabla(\mu)) \cong 0.
$$
Futhermore, if for all $i$ and $j > i$ neither $\lambda_i - \lambda_j -i +j$ nor
$\mu_i - \mu_j -i +j$ is divisible by $l$, and $\lambda$ and $\mu$ are
in the same block of $G$, then
$$
\Ext^{d(\lambda) - d(\mu)}_{S(n,r)}
(L(\lambda), \nabla(\mu)) \cong 
\Ext^{d(\lambda) - d(\mu)}_{S(n,r)}
(\nabla(\lambda), \nabla(\mu)) \cong k.
$$
\end{thm}

\begin{rem}
In this paper we have tried to avoid using alcove-geometric notation,
for simplicity of exposition.  However, the conditions in the second
part of this theorem
(and Theorem \ref{parker}) can be restated as requiring $\lambda$ (and
$\mu$) to lie in the interior of an alcove, with $\lambda$ and $\mu$
in the same orbit of the affine Weyl group $W_l$ (under the usual \lq
dot\rq\ action on weights). When thus stated both results are valid
for an arbitrary reductive algebraic group.
\end{rem}

If $l-3 \ge d(\lambda) - d(\mu)$ then an application of
Corollary \ref{cor:sameext} would give the corresponding result for the
Hecke algebra --- but this is only going to be true for partitions
that are close together. 

The following very nice result of Wen \cite[theorem 8.3, lemma
  8.3.3]{wen} 
gives the value of the $\Ext$ groups (in small
degrees) between
induced modules when we have weights which are close together and
related by a particular type of element in the affine Weyl
group. Recall that a weight $\lambda$ is \emph{strictly dominant} if
$\lambda-(l-1)\rho$ is dominant. 
(Again, this result was stated and proved for the classical case, but
the quantum version follows in the same way.)

\begin{thm}[Wen]\label{wen}
Let $\lambda \in \partn(n,r)$ be inside an alcove and strictly
dominant.
Let $S$ be a commutative subset of the reflections which fix the walls
in the closure of the alcove containing $\lambda$.
Let $w$ be in the subgroup of $W_l$ generated by $S$ 
such that $w\cdot\lambda < \lambda$. 
%
Set $d = d(\lambda) - d(w\cdot\lambda)$.
We have
for each $i$ with $0 \le i \le d$ 
that
$$
\dim\Ext^i_{G}(L(\lambda), \nabla(w \cdot \lambda)) 
= \left\{\begin{array}{ll}
1 &{\mbox{if $d = i$}}\\
0 &{\mbox{otherwise.}}
\end{array}
\right.
$$
$$
\dim\Ext^i_{G}(\nabla(\lambda), \nabla(w \cdot \lambda)) 
= \binom{d}{i}.
$$
\end{thm}
This result gives a special case of Theorem \ref{ryom} but also finds
all the $\Ext$ groups for weights satisifying the special conditions
above.  The usefulness of the result is limited by the size of the
subset $S$ and this has an upper bound of $\frac{n}{2}$. As usual, we
have only stated the type $A$ version; the general result holds for
any reductive group.

In a similar spirit to Theorem \ref{ryom} we have \cite[Corollaries 4.4 and 4.5]{par2}
 
\begin{thm}\label{parker}
Suppose $\lambda, \mu \in \partn$. If $m \ge d(\lambda) + d(\mu)$
then we have
$$
\Ext^{m}_{S(n,r)}
(L(\lambda), L(\mu)) \cong 
\Ext^{m}_{S(n,r)}
(\nabla(\lambda), \Delta(\mu))
\cong 
0.
$$
Furthermore, if $\lambda_i - \lambda_j -i +j$ is not divisible by
$l$ for all
$i$ and $j > i$ and if
$\mu_i - \mu_j -i +j$ is not divisible by $l$ for all
$i$ and $j > i$ then  
we have 
$$
\Ext^{m}_{S(n,r)}
(L(\lambda), L(\mu)) \cong 
\Ext^{m}_{S(n,r)}
(\nabla(\lambda), \Delta(\mu))
\cong 
k.
$$
\end{thm}

Using these results, the second author has determined \cite[Theorems
  5.8 and 5.9]{par2} the global dimension of the Schur algebra
$S(n,r)$ when $l>n$ or $l=n$ and $r\equiv 0 \pmod l$: 

\begin{thm}
If $l>n$ then the global dimension of $S(n,r)$ is
$$2(n-1)\left\lfloor \frac{r}{l}\right\rfloor.$$
The global dimension of $S(l,ml)$ is
$$2(l-1)m.$$
\end{thm}

When $r\leq n$ the global dimension of $S(n,r)$ has been calculated by
Totaro \cite{tot} in the classical case, and Donkin \cite[Section
  4.8]{donkbk} in the quantum case.

Another sometimes useful fact about $\Ext$ groups is the following.
\begin{lem}\label{extsum}
If $\mu \ne \lambda$ then
$$\sum (-1)^{i} \dim \Ext^i_G(\nabla(\mu), \nabla(\lambda)) = 0$$
\end{lem}
This is a consequence of \cite[II, 6.21(6)]{jantzen}.  We use this to
show that there are two dimensional $\Ext^1$ groups for $G$ for
$l\ge3$ and $n\ge 3$, even for weights which are close together.
Suppose $\lambda=(\lambda_1,\lambda_2,\lambda_3)$ is a dominant weight
for $G$ with $l\ge3$.  We further suppose that
$\lambda_1-\lambda_2 \not \equiv -1\pmod l$, $\lambda_2-\lambda_3 \not
\equiv -1\pmod l$, and $\lambda_1-\lambda_3 \not \equiv -2\pmod l$,
(i.e. that $\lambda$ lies in the interior of an alcove), and that
$d(\lambda) \ge 2$.  Consider a dominant weight $\mu < \lambda $ which
is in the same block as $\lambda$, does not satisfy $\lambda -\mu \in
\N \alpha$ where $\alpha$ is a root, and satisfies $d(\lambda) -
d(\mu) =2$. (Such a weight exists as $d(\lambda) \ge 2$.)  We then
have that $\Hom_{G}(\nabla(\lambda), \nabla(\mu))\cong k$,
using the results of \cite{coxpar} or Theorem \ref{kopthm},
and that $\Ext^2_{G}(\nabla(\lambda), \nabla(\mu))\cong k$ and
$\Ext^i_{G}(\nabla(\lambda), \nabla(\mu))\cong 0$ for $i \ge3$
using Theorem \ref{ryom}. Thus using Lemma \ref{extsum} we have that
$$\Ext^1_{G}(\nabla(\lambda), \nabla(\mu))\cong k^2.$$
 So there
are two-dimensional $\Ext^1$ groups for $l\ge3$, even for
weights which are close together. 

The above example shows that calculating $\Ext$ groups for modules for
groups of rank at least three is quite complicated in general. Indeed,
this problem may be intractable for groups of large rank. There is
still possibly some hope of calculating $\Ext^1$ between induced modules
for $\GL_3(k)$, but the dimension of the $\Ext^1$ groups may be
unbounded for $n=3$.  This is in contrast to the $n=2$ case where the
$\Ext^1$ groups between induced modules have dimension at most one.
Theorem \ref{wen} shows that the dimension of the $\Ext^1$ group can
be as large as we like, as long as we take $n$ large enough.

\section{Ext-groups involving simple modules}\label{simples}
In the paper \cite{klsh} Kleshchev and Sheth ask when
$$\Ext^1_H(D^\lambda, D^\mu) \cong \Hom_H( \rad S^\lambda, D^\mu)?$$
One sufficient condition for the equality to
hold will be if $\Ext_H^1(S^\lambda, D^\mu) \cong 0$. 
Kleshchev and Sheth prove the following for the classical case.
\begin{thm}[{\cite[theorems 2.9 and 2.10]{klsh}}]
Suppose $p\ne 2$,  $\lambda \not >\mu$ and
$\lambda$ and $\mu \in \rpartn(p-1,r)$ then
$$\Ext_H^1(S^\lambda, D^\mu) \cong 
\left\{
\begin{array}{ll}
k &\mbox{if $\lambda = \mu = \epsilon_r$ and $r \ge p$}\\
0 &\mbox{otherwise.}
\end{array}
\right.
$$ where $\epsilon_r$ is the Mullineax map applied to the partition
$(r,0,0,\ldots)$ and
$$\Ext^1_H(D^\lambda, D^\mu) \cong \Hom_H( \rad S^\lambda, D^\mu).$$
\end{thm}

This motivates the following lemma.
\begin{lem}
If $\lambda$, $\mu \in\rpartn(n,r)$ and $l\ge 4$ then
$\Ext^1_H(S^\lambda, D^\mu) \not\cong 0$ implies that $m(\mu)' <
\lambda$.
\end{lem}
\begin{proof}
We will apply the Schur functor to the following short exact sequence
for the tilting module $T(\lambda)$:
$$ 0 \to K \to T(\lambda) \to \nabla(\lambda) \to 0.$$
By \cite[4.4(2)]{donkbk}
we get the following short exact sequence
$$
0 \to fK \to Y^{\lambda'} \otimes \sgn \to S^\lambda  \to 0$$
where $fK$ has a filtration
by Specht modules $S^{\nu_i}$ with $\nu_i < \lambda$.

Now
$$
\Ext_H^1(Y^{\lambda'} \otimes \sgn, D^\mu) \cong
\Ext_H^1(Y^{\lambda'}, D^\mu\otimes \sgn) \cong 0
$$ 
as $Y^{\lambda'}$ is projective if $\lambda \in \rpartn(n,r)$.
Thus we get the following sequence
\begin{align*}
0 \to
\Hom_H(S^\lambda, D^\mu) \to
\Hom_H(Y^{\lambda'} \otimes \sgn, D^\mu) \to
\Hom_H(fK, D^\mu) \to
\Ext_H^1(S^\lambda, D^\mu) 
\to 0
\end{align*}
Thus $\Ext_H^1(S^\lambda, D^\mu) \not \cong 0$ implies that
$\Hom_H(fK,D^\mu)$ is non-zero. 
Now $$\Hom_H(fK,D^\mu) \cong \Hom_S(gfK, L(m(\mu)'))$$ using Theorem
\ref{thm:donk}.
But $l \ge 4$ and $fK$ has a Specht filtration with $(fK:S^\nu)=
(K:\nabla(\nu))$ and so $gfK \cong K$.
Thus 
$\Hom_S(gfK, L(m(\mu)'))$ can only be non-zero if $L(m(\mu)')$ is in the
head of one of the $\nabla(\nu)$ (with $\nu < \lambda$)
appearing in a filtration of $T(\lambda)$.
Let $\nabla(\nu)$ be such a $\nu$ with head $L(m(\mu)')$.
Thus $m(\mu)' \le \nu <\lambda$.
\end{proof}

We can thus generalise the result of Kleschev and Sheth to the
following
\begin{cor}\label{whencong} If $\lambda$, $\mu \in \rpartn(r)$, $l\ge4$  and
  $m(\mu)' \not < \lambda$ then 
$$\Ext^1_H(D^\lambda, D^\mu) \cong \Hom_H( \rad S^\lambda, D^\mu).$$
\end{cor}

Kleshchev and Martin  have conjectured that for $l\geq 3$ we have
$$\Ext^1_{k\Sym_r}(D^{\lambda}, D^{\lambda}) = 0.$$ 
Unfortunately Corollary \ref{whencong} does not help us show that 
$\Ext^1_H(D^{\lambda}, D^{\lambda}) = 0$, as
$m(\lambda)' \le \lambda$. This can be seen by considering the Specht
module $S^{m(\lambda)'}$. We claim that $D^{\lambda}$ occurs once as a
composition factor (indeed it is the socle of $S^{m(\lambda)'}$) and
so $\lambda \ge m(\lambda)'$.
Now $[S^{m(\lambda)'}: D^{\lambda}]
= [S^{m(\lambda)'}\otimes \sgn: D^{\lambda}\otimes \sgn]
= [(S^{m(\lambda)})^*: D^{m(\lambda)}]
= [S^{m(\lambda)}: D^{m(\lambda)}]
=1 $.

We can now deduce that the condition $\lambda \not > m(\mu)'$ implies
$\lambda \not > \mu$ but that the reverse implication does not hold in
general.

Hemmer \cite{hemext, hemrow} also gets results concerning
extensions between simples. In particular, he determines the non-split
extensions between simple modules when the simples involved are
completely splittable.


\bibliographystyle{amsalpha}

\providecommand{\bysame}{\leavevmode\hbox to3em{\hrulefill}\thinspace}
\providecommand{\MR}{\relax\ifhmode\unskip\space\fi MR }
\providecommand{\MRhref}[2]{%
  \href{http://www.ams.org/mathscinet-getitem?mr=#1}{#2}
}
\providecommand{\href}[2]{#2}

\end{document}